\documentclass[a4paper,twocolumn]{article} 
\pdfoutput=1
\usepackage{amsmath}
\usepackage{authblk}
\usepackage{booktabs}
\usepackage{color}
\usepackage{graphicx}
\usepackage{listings}
\usepackage{subfig}
\usepackage{url}
\usepackage[pdfborder={0 0 0},colorlinks=true,citecolor=red,linkcolor=blue]{hyperref}

\definecolor{light-gray}{gray}{0.95}
\title{Towards a High Fidelity Direct Transcription Method for Optimisation of Low-Thrust Trajectories}
\author{Chit Hong Yam\footnote{\href{mailto:chit.hong.yam@esa.int}{chit.hong.yam@esa.int}}}
\author{Dario Izzo\footnote{\href{mailto:dario.izzo@esa.int}{dario.izzo@esa.int}}}
\author{Francesco Biscani\footnote{\href{mailto::bluescarni@gmail.com}{bluescarni@gmail.com}}}
\affil{European Space Agency -- Advanced Concepts Team \\ European Space Research and Technology Centre (ESTEC)}

\begin{document}

\maketitle

\begin{abstract}
We build upon some new ideas in direct transcription methods developed within the Advanced Concepts Team
to introduce two improvements to the Sims-Flanagan transcription for low-thrust trajectories. The obtained
new algorithm is able to produce an operational trajectory accounting for the real spacecraft dynamics and
adapting the segment duration on-line improving the final trajectory optimality.
\end{abstract}

\section{Introduction}
A direct optimisation method proposed by Sims and Flanagan \cite{GOlowthrust:Sims99} suggests that low-thrust
trajectories can be modeled as a series of impulsive $\Delta$V connected by conic arcs. The method is fast and
robust and has been applied in previous works for preliminary mission design \cite{JIMO:Yam1,JIMO:Yam2}. However with its impulsive
$\Delta$V transcription, the Sims-Flanagan method can fail to accurately represent the actual dynamical
model unless the number of impulses is increased and thus at the cost of slowing down the overall optimisation. 

In this paper, we study new transcription methods to improve the accuracy of the Sims and Flanagan model without
increasing the dimension of the problem. The works extends some of the ideas presented at the fifth
international meeting on celestial mechanics  CELMEC V \cite{celmec}. The first improvement is to replace
the impulses with continuous thrust, where low-thrust arcs are numerically propagated. The magnitude and
direction of the thrust are part of the optimisation variables and are assumed to be constant throughout a
segment. Perturbations can be included, in the propagation, to further improve the fidelity of the model.
The modification introduces a performance penalty due to the higher computational costs of the integration
with respect to a simple Keplerian propagation between impulses. In order to tackle this issue we introduce
the use of Taylor integration \cite{taylorint} methods in place of the commonly employed Runge-Kutta-Fehlberg
scheme, reducing the performance loss by almost one order of magnitude.

A second improvement we introduce is to allow the time mesh to be optimised together with the trajectory.
While this is a long unsolved issue in direct method for trajectory optimisation, we manage to obtain an
efficient algorithm by introducing the Sundman transformation\cite{sundman19112}, in which the independent
variable is changed from time to $s$, where the time variation of $s$ is inversely proportional to the radial
distance. By doing so, and adding only one constraint to the optimisation problem, segments are automatically
distributed more densely near the central body (where speed is usually higher)
along the optimal solution and thus on-line mesh adaptation is obtained at the cost of an acceptable
performance loss. We present a numerical example to compare results between the original and the new methods.
The resulting tool has the further advantage of being suitable for different phases of the mission
design, from preliminary, where global optimisation methods need a rather simple and low-dimensional
transcription, to operational where dynamics need to be accounted for in a precise manner and optimality is sought.
 
\begin{figure}[htp]
\centerline{\includegraphics[width=8cm]{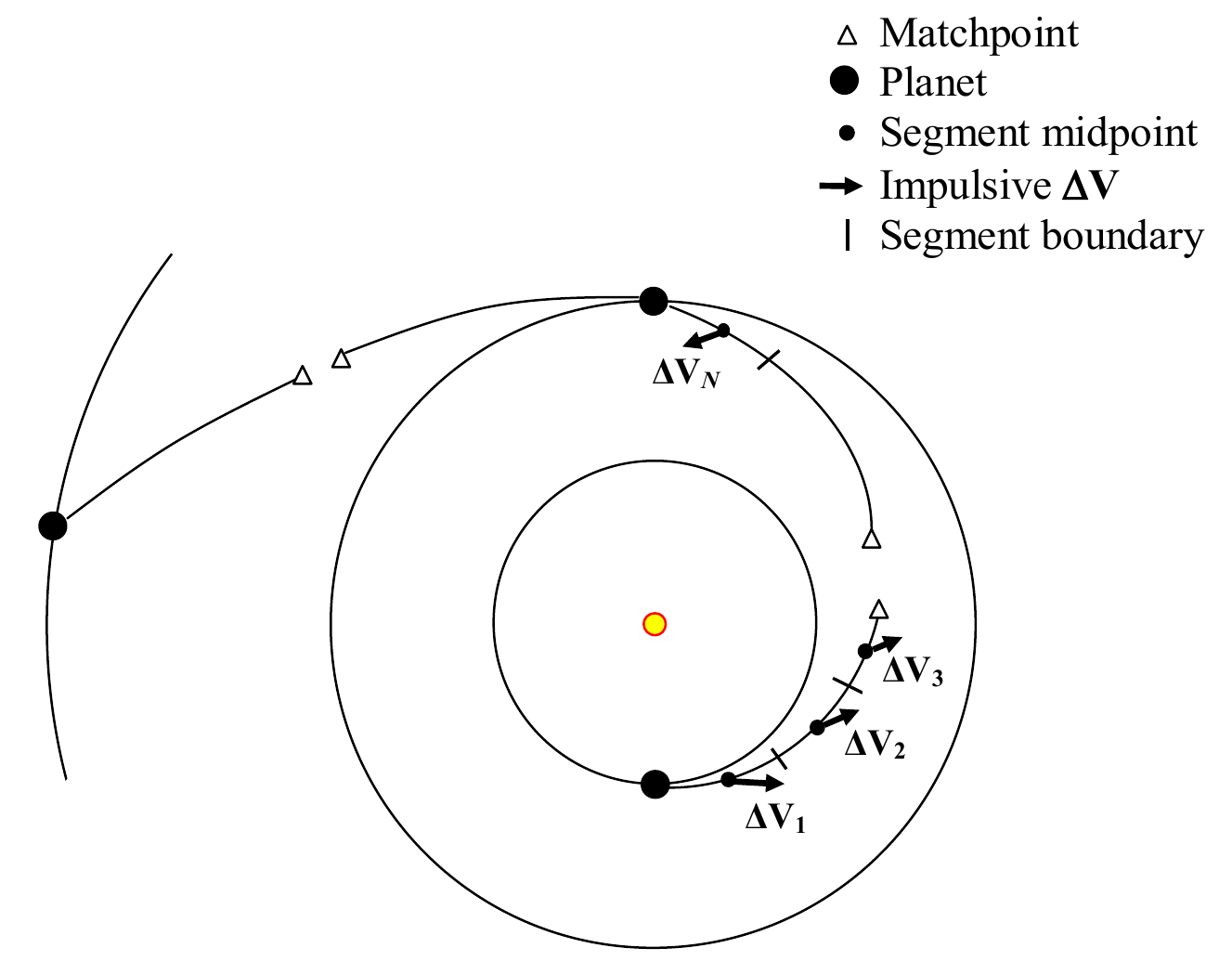}}
\caption{\textnormal{Impulsive $\Delta$V transcription of a low-thrust trajectory, after Sims and Flanagan 
\cite{GOlowthrust:Sims99}}}
\label{fig_sim}
\end{figure}

\section{The Original Sims-Flanagan Model}

 In 1999, Sims and Flanagan proposed a direct method for optimising low-thrust trajectories, where later the software packages GALLOP \cite{McConaghy03} and MALTO \cite{Sims06} are developed based on this model. Figure~\ref{fig_sim} briefly illustrates such a trajectory model. The whole trajectory is divided into legs which begin and end with a planet. Low-thrust arcs on each leg are modeled as sequences of impulsive maneuvers $\Delta$V, connected by conic arcs.  We denote the number of impulses (which is the same as the number of segments) with $N$. The $\Delta$V at each segment of equal duration should not exceed a maximum magnitude, $\Delta$V$_{max}$, where $\Delta$V$_{max}$ is the velocity change accumulated by the spacecraft when it is operated at full thrust during that segment:

\begin{equation}
\Delta V_{max} = ( F_{max}/m) (T_f - T_0) / N
\label{eq:eq1}
\end{equation}
  
where $F_{max}$ is the maximum thrust of the low-thrust engine, $m$ is the mass of the spacecraft,  $T_0$ and $T_f$ is the initial and final time of a leg. The spacecraft mass is propagated using the rocket equation \cite{rocketeq}:
\begin{equation}
    m_{i+1} = m_i \text{ exp}(-\Delta V_i/g_0I_{sp})
\label{eq:rocket}
\end{equation}     
where the subscript $i$ denotes the mass and $\Delta$V on the $i$-th segment, $g_0$ is the standard gravity (9.80665 m/s$^2$), and $I_{sp}$ is the specific impulse of the low-thrust engine.

	At each leg, trajectory is propagated (with a two-body model) forward and backward to a matchpoint (usually halfway through a leg), where the spacecraft state vector becomes \textbf{S}$_{mf}$ = $\left\{ r_x, r_y, r_z, v_x, v_y, v_z, m \right\}$$_{mf}$ (and similarly for \textbf{S}$_{mb}$), where $r$ and $v$ are respectively the position and velocity of the spacecraft and the subscripts represents the Cartesian $x$, $y$, $z$ components.  The forward- and backward-propagated half-legs should meet at the matchpoint, or the mismatch in position, velocity, and mass: 

\begin{equation}
     \textbf{S}_{mf} - \textbf{S}_{mb}
          = \left\{\Delta r_x,\Delta r_y,\Delta r_z,\Delta v_x,\Delta v_y,\Delta v_z, \Delta m \right\}     
\label{eq:eq2}
\end{equation}                
should be less than a tolerance in order to have a feasible trajectory. 

The problem is transcripted into a nonlinear programming problem (NLP), where the objective is to maximize the
final spacecraft mass subjected to the constraints on the maximum $\Delta$V and the state mismatch, while the
decision variables of the problem are listed below:
\begin{itemize}
  \item the departure epoch $T_0$
  \item the departure velocity relative to the earth $V_{\infty}$
  \item for each leg and each segment, the magnitude of the impulse
    and direction
  \item for each swingby, the incoming and outgoing velocities
    relative to the planet
  \item for each swingby $j$, the swingby epoch $T_j$
  \item the arrival epoch $T_f$
\end{itemize}

For a rendezvous mission, the arrival velocity to the destination is not included in the set of variables,
as it is, by construction of the model, zero relative to the planet. To solve the NLP, we use a software
package called SNOPT \cite{SNOPT02}  \cite{SNOPT06}, which implements sequential quadratic programming (SQP).

\begin{figure}[ht]
\centering
\includegraphics[width=8cm]{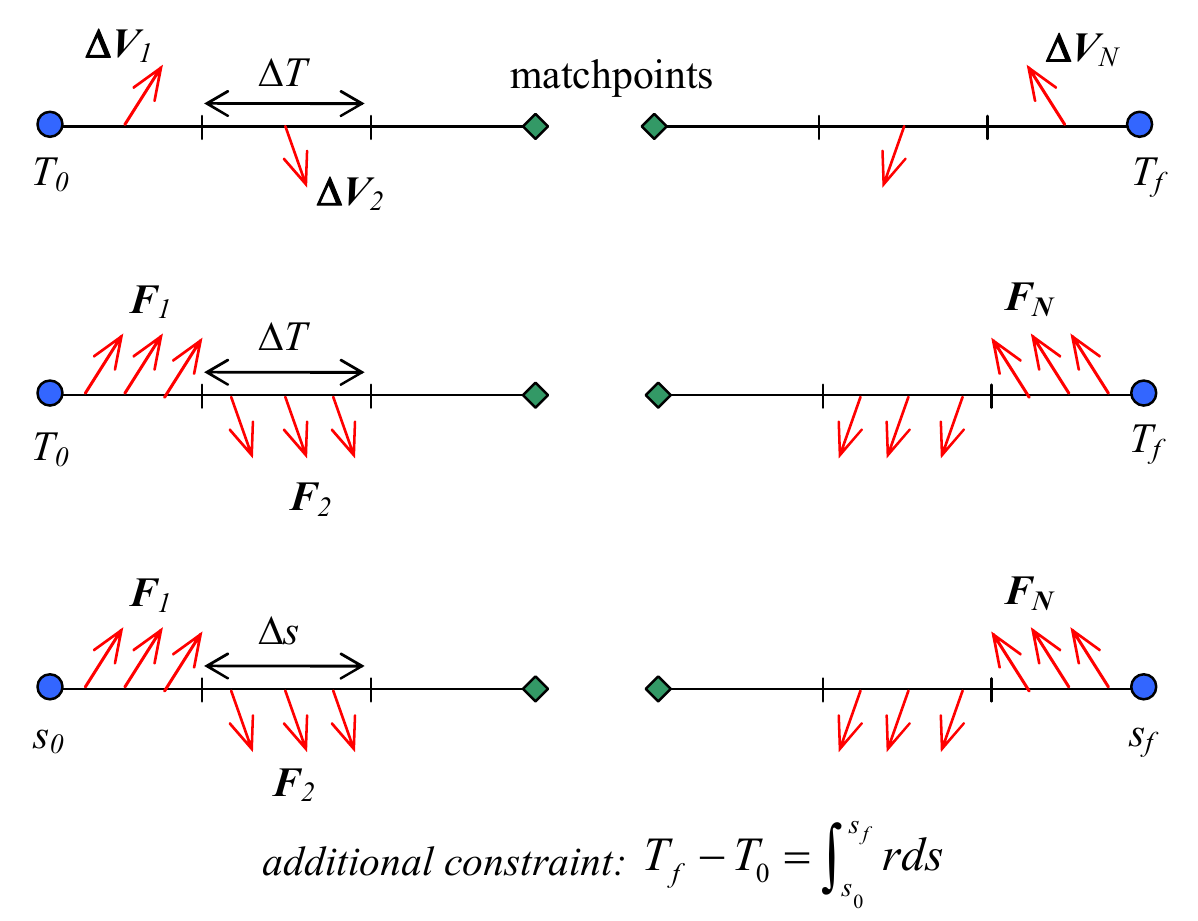}
\caption{Schematics of the three trajectory models. Top: Impulisve $\Delta$V; Middle: Continuous thrust is
the time-space; Bottom: Continuous thrust in the s-space} \label{fig:compare3models}
\end{figure}

\section{Improvement to Trajectory Model}
From the beginning, the use of the Sims-Flanagan model is limited to preliminary mission design, in which
the results are not expected to be accurate up to the operation level. In terms of the fidelity of the model,
there are two areas in the Sim-Flanagan model that can lead to loss in accuracy: (1) the use of impulses; and
(2) the insufficient number of segments.  To address these issues, we introduce two improvements:


\begin{itemize}
\item Impulsive $\Delta$V are replaced by continuous thrust to improve the fidelity of the trajectory dynamics. In order
to keep a reasonably low computing time, we employ a Taylor integration scheme showing an order of magnitude performance 
gain with respect to the classical Runge-Kutta methods.
\item An adaptive time mesh is obtained on the segments via the Sundman transformation to improve on the optimality of the
 final trajectory. 
\end{itemize}
\subsection{The Taylor Integration Method}

When replacing the original impulsive $\Delta V$ transcritpion with a continuous fixed thrust transcription,
the optimization process becomes slower as the optimization relies on a numerical integration scheme of a
higher complexity with respect to a simpler ballistic arc solver. Efficiency is essential
to keep the CPU time penalty at a minimum level. For this purpose, we report the comparison,
in terms of CPU time and accuracy, among classical Runge-Kutta-Fhelberg methods and the Taylor integration \cite{taylorint}. The tests have been done having in mind the typical algorithm call done during an 
optimization procedure that uses our approach. An improvment of one order of magnitude in CPU time (while 
keeping the integration accuracy to the same level) is found.

\paragraph{Comparison set-up}
We consider the set of differential equations describing the motion of a spacecraft subject to a fixed thrust 
force in the interplanetary medium. This \lq fixed thrust problem\rq\ is at the basis of the ideas on direct transcription 
methods presented by some of these authors during the  fifth international meeting on celestial mechanics 
CELMEC V (\cite{celmec}) and that motivated the current paper. Since in the proposed new direct transcription 
method the fixed thrust problem needs to be solved a large amount of times (in each segment) and with diverse 
initial conditions and thrust vectors we focus the algorithmic comparison to those cases representative of such 
a process. The equations, in a non dimensional form, are the following:

\begin{equation}
\begin{array}{l}
r  =  \sqrt{x_1^2 + x_2^2+x_3^2} \\
\dot x_1 = x_4\\
\dot x_2 = x_5\\
\dot x_3 = x_6\\
\dot x_4 = - x_1 / r^3 + u_1 \\
\dot x_5 = - x_2 / r^3 + u_2 \\
\dot x_6 = - x_3 / r^3 + u_3 
\end{array}
\label{eq:motion}
\end{equation}

The mass is not considered for the purpose of this comparison, but will be included in the trajectory model 
described later.
To test the integration schemes, $N=10000$ different Cauchy problems have been generated at random considering 
$x_i(0), u_i$ uniformly distributed in $x_i(0) \in [0.1,2]$ and $u_i \in [0.0001,0.01]$. The final integration 
time has been also set to be random and $t_f \in [\pi / 20,10\pi]$. The same problems where solved using a 
Runge-Kutta-Fhelberg integration scheme (in the implementation of the GAL libraries \cite{gal}) and a Taylor 
integration scheme (implemented using the tool \lq\lq taylor\rq\rq\ \cite{taylorweb}). In order to test the speed
and the precision of 
the solvers, we propagate each problem from $x_i(0)$ for $t_f$, we then take the result and propagate backwards 
for $t_f$ reaching the point $x_i^f$. By doing this, as we know the exact result of the propagation that is $x_i(0)$, 
we evaluate the precision of the propagation defining the
propagation error as $err = \sum_{i=1}^6 \left(x_i(0)-x_i^f \right)^2$. 
Each algorithm is tested on the same set of randomly generated Cauchy problems. In all cases, no minimum step 
size is used and the same parameter $\epsilon$ is passed to the RKF integrators as the absolute error, and to 
the Taylor integrator as both absolute and relative error. The initial trial stepsize of 0.1 is set to the RKF 
integrators.

\paragraph{Results}

From the results oulined in Table \ref{tab:comp} it is clear that the Taylor integrator is outperforming the RKF 
both in speed and accuracy confirming in the low-thrust fixed direction problem the same performance gain levels already 
reported in past literature \cite{taylorint,scott2010high}. From the table we may also, empirically, establish that
$\epsilon = 10^{-10}$ is a good compromise between speed and 
accuracy and can thus be used as a default parameter to call the Taylor integrator. The speed gained by employing 
the Taylor integrator is roughly one order of magnitude.

\begin{table*}[ht]
\caption{Algorithm Performance. Speed is measured in seconds and refers to all the $N$ integrations (forward+backward). Max. Err. is the maximum integration error made in the $N$ propagations. Note that this is a real error, not an estimation as explained above. }
\begin{center}
\begin{tabular}{c|cc|cc|cc}
 $\epsilon$ & \multicolumn{2}{c}{RKF 5(6)} & \multicolumn{2}{c}{RKF 7(8)} & \multicolumn{2}{c}{Taylor}\\\hline
 & Speed (s) & Max. Err. &  Speed (s) & Max. Err. &  Speed (s) & Max. Err. \\
 1.00e-03&6.10e-01&1.45e+01&8.73e-01&4.69e+01&4.59e-01&1.29e+01\\
1.00e-04&8.07e-01&3.12e+01&1.13e+00&7.51e-01&6.59e-01&9.71e-04\\
1.00e-05&1.01e+00&2.90e+00&1.40e+00&2.39e-01&8.19e-01&4.66e-04\\
1.00e-06&1.30e+00&8.08e-01&1.63e+00&2.21e-02&9.33e-01&1.82e-07\\
1.00e-07&1.81e+00&6.24e-02&1.93e+00&6.10e-03&1.12e+00&1.21e-07\\
1.00e-08&2.61e+00&1.39e-03&2.37e+00&1.18e-03&1.36e+00&4.75e-11\\
1.00e-09&3.69e+00&7.06e-05&3.03e+00&1.49e-04&1.53e+00&4.90e-11\\
1.00e-10&5.61e+00&1.48e-05&3.94e+00&1.32e-05&1.98e+00&6.39e-15\\
1.00e-11&8.76e+00&1.67e-06&4.95e+00&2.68e-07&2.28e+00&2.15e-18\\
1.00e-12&1.34e+01&1.73e-07&6.14e+00&2.20e-08&2.45e+00&4.78e-18\\
1.00e-13&2.06e+01&1.87e-08&8.24e+00&2.78e-09&2.65e+00&9.97e-19\\
1.00e-14&3.25e+01&1.96e-09&1.09e+01&6.11e-10&2.91e+00&3.59e-20\\
1.00e-15&5.41e+01&7.35e-10&1.52e+01&1.76e-09&3.17e+00&1.97e-19\\
1.00e-16&1.30e+02&2.58e-09&2.68e+01&6.77e-10&3.80e+00&1.61e-19
 \end{tabular}
\end{center}
\label{tab:comp}
\end{table*}%

\subsection{The Sundman Transformation}

In his celebrated paper \cite{sundman19112} Karl Sundman introduces a simple differential transformation for the time
variable, to regularize the otherwise singular three body problem. The Sundman transformation dilates the time 
metric introducing a new variable $s$ defined through the relation $ds = dt / r$ guaranteeing an asymptotically 
slower flow near the singularities. In a trajectory propogation this same property turns out to be quite useful 
if equally spaced segments are considered in the $s$ domain rather than in the $t$ domain. Let us for example 
consider Eq.(\ref{eq:motion}) and use the Sundman transformation, we obtain the following set of equations:

$$
\begin{array}{l}
r =  \sqrt{x_1^2 + x_2^2+x_3^2} \\
\dot x_1 = x_4 r\\
\dot x_2 = x_5 r\\
\dot x_3 = x_6 r \\
\dot x_4 = - x_1 / r^2 + u_1 r \\
\dot x_5 = - x_2 / r^2 + u_2 r \\
\dot x_6 = - x_3 / r^2 + u_3 r \\
\dot t =  r
\end{array}
$$

\begin{figure*}[ht]
\centering
a) \includegraphics[width=7cm]{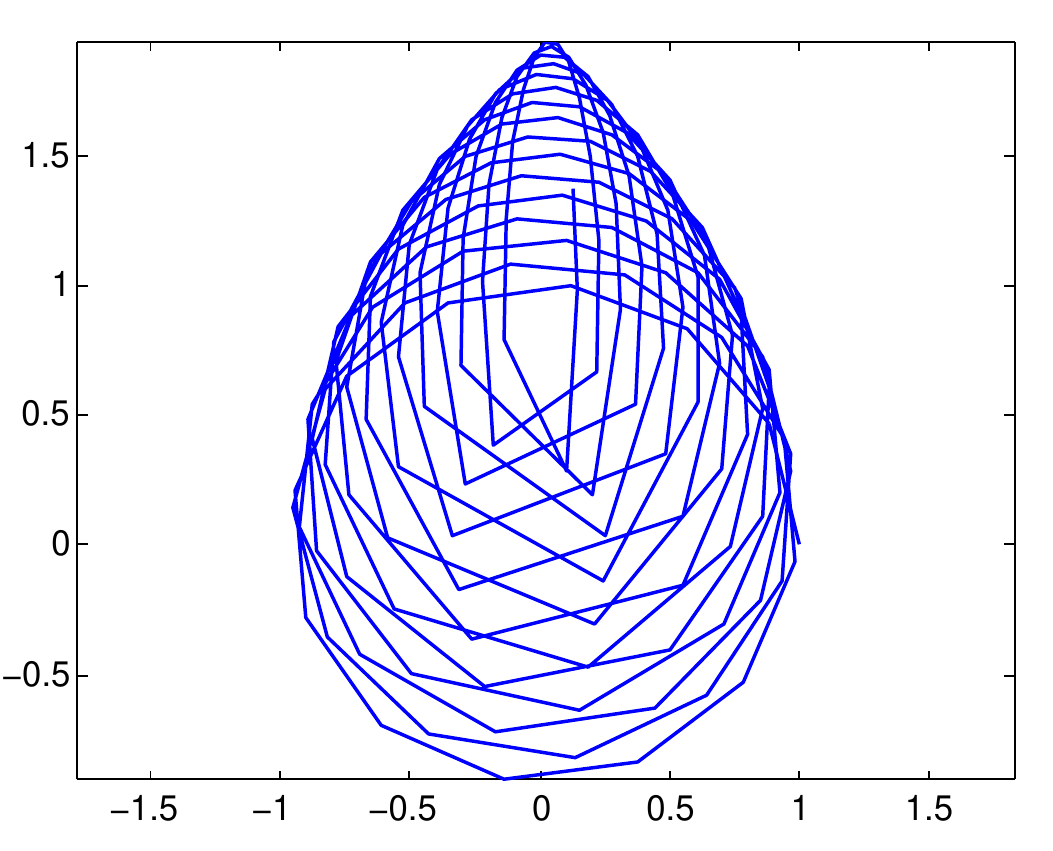}
b) \includegraphics[width=7cm]{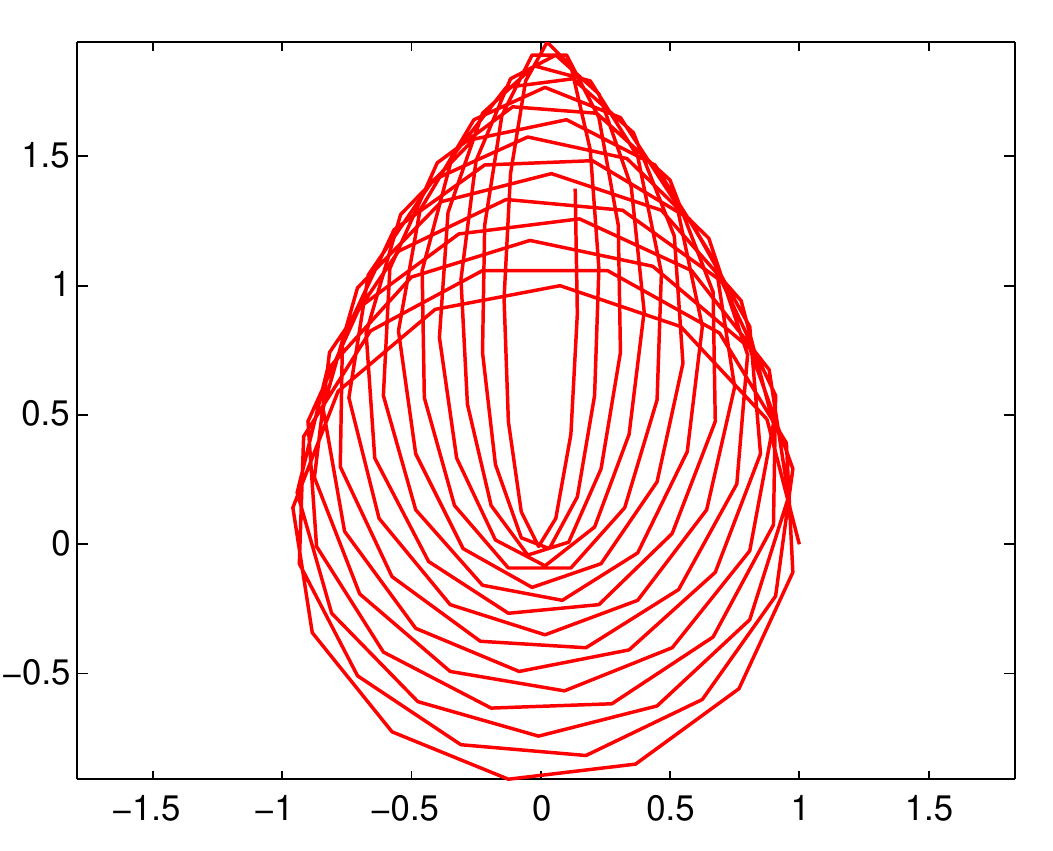}
\caption{A trajectory sampled with the same number of equally a) time spaced segments b) s-spaced segments} 
\label{fig:time1}
\end{figure*}

\begin{figure*}[ht]
\centering
a) \includegraphics[width=7cm]{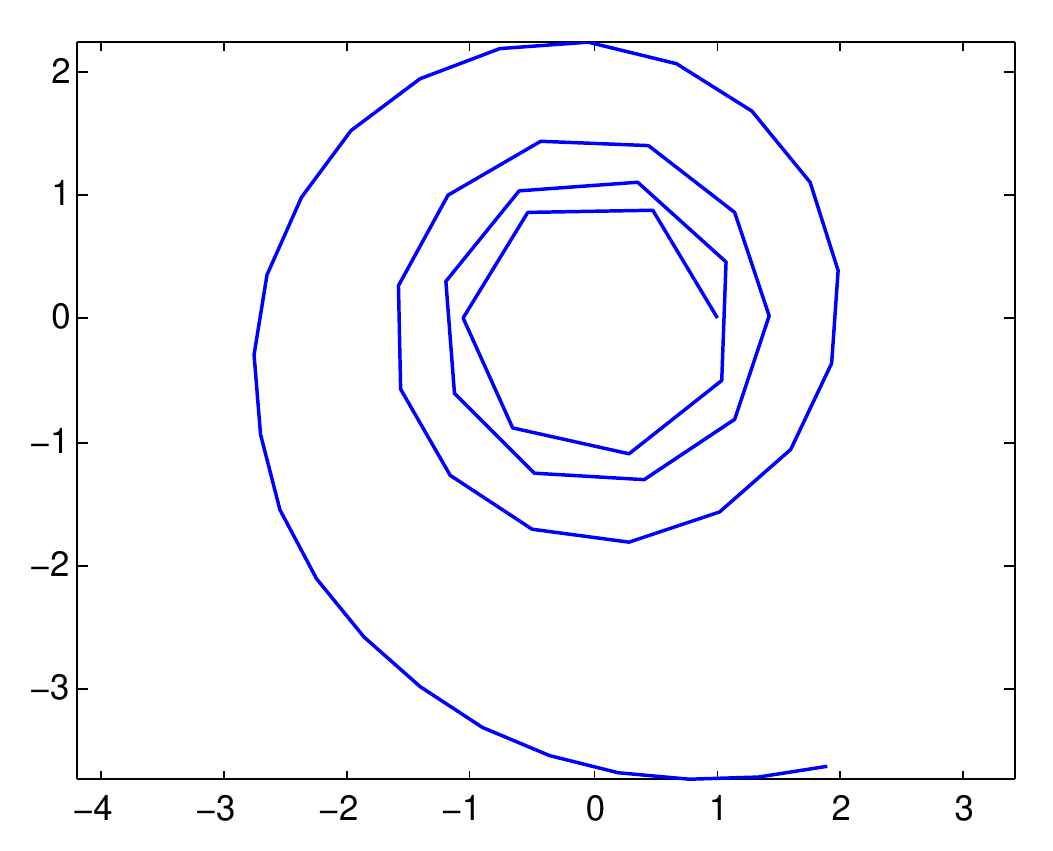}
b) \includegraphics[width=7cm]{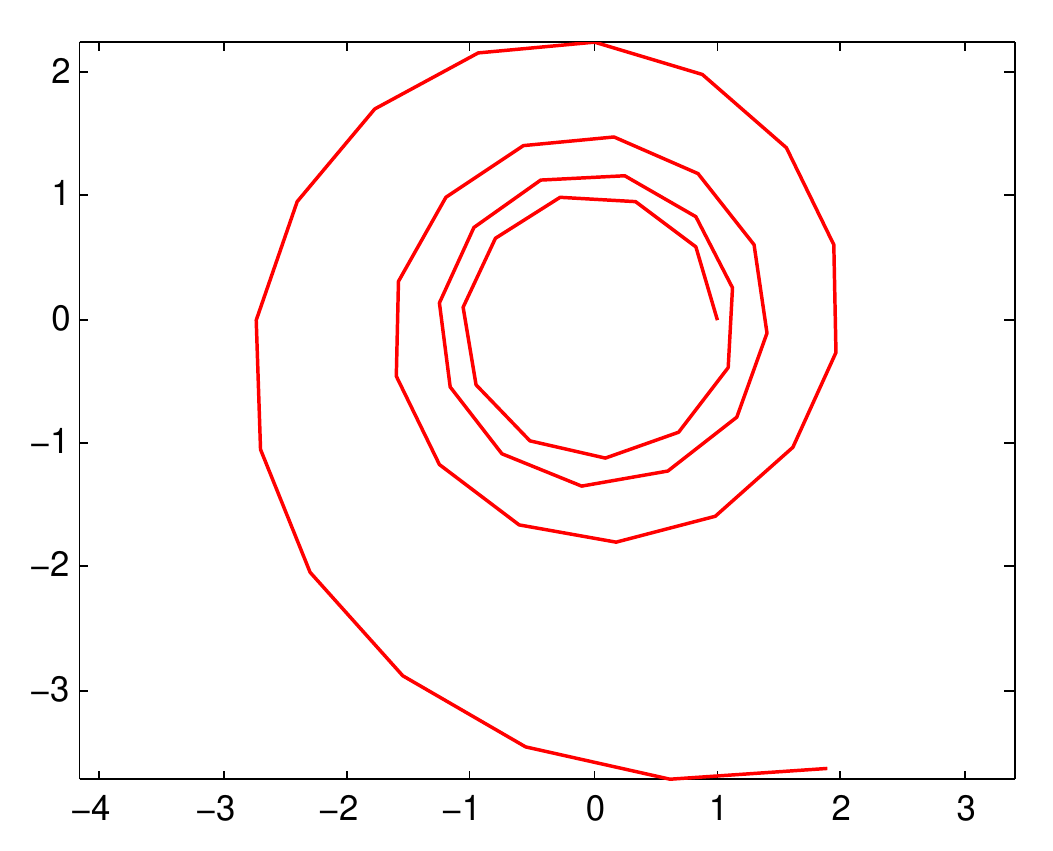}
\caption{A trajectory sampled with the same number of equally a) time spaced segments b) s-spaced segments} 
\label{fig:time2}
\end{figure*}

\begin{table*}[ht]

\begin{minipage}[t]{8.5cm}
\subfloat[][Parameters for an Earth-Mercury Rendezvous Mission.]
{
\centering
\begin{tabular}{lr}
\toprule
\multicolumn{1}{l}{\!Parameters\!}
& \multicolumn{1}{r}{Values} \\  
\midrule

 Initial mass of the spacecraft		&    $660.0$ kg  \\
 Maximum thrust    			&    $92.3$ mN \\
 Specific impulse  			&    $3,337$ s  \\
 Launch date  			&    Apr. 9, 2007  \\
 Arrival date  			&    Aug. 22, 2013  \\
 Launch $V_\infty$ 	&    $\leq 2.0$ km/s    \\
 Time of flight     & 6.37 years   \\
\bottomrule
\end{tabular}
\label{tab-NEPparam}
}
\end{minipage}
\begin{minipage}[t]{8.5cm}
\subfloat[][Optimal Final Mass for the Earth-Mercury Mission. Note that only the continuous trajecctories are feasible
and thus can be used up to late design phases.d]
{
\centering
\begin{tabular}{lcc}
\toprule
\multicolumn{1}{l}{\!Method\!}
& \multicolumn{1}{c}{No. of}
& \multicolumn{1}{c}{Optimal Final} \\
& Variables & Mass, kg \\

\midrule
Impulsive                    & $97$ &  $392.9$ \\
Continuous \textit{t}-space        & $97$  & $382.4$ \\
Continuous \textit{s}-space  & $98$ & $387.0$ \\
\bottomrule
\end{tabular}
\label{tab-EMercuryMf}
}
\end{minipage}
\caption{Trajectories parameters and results.}
\end{table*}

To demonstrate the effect of such a transformation on a numerical mesh, we take a circular orbit of radius one and
propagate it forward with a constant thrust aligned along the $x$ axis. In Figure \ref{fig:time1}a we visualize the
obtained orbit using a uniform sampling in time, while in Figure \ref{fig:time2}b the same trajectory is visulaized
using the same number of samples, but equally spaced in the Sundman variable domain. The same is done in Figures
\ref{fig:time2}a-b for a constant tangential thrust. These pictures clearly show the problem with using points
equally spaced in time to define a mesh for a numerical algorithm: due to the conservation of energy the closer we get
to the singularity the more potential energy we lose and thus acquire in terms in kinetic energy. This pumps up the body 
velocity substantially creating an unequal distribution of segments length bound to create numercial difficulties. The use
of the $s$ variable is one of the possible transformations able to alleviate such a problem.

\section{New trajectory models}
The implementation of the ideas reported above leads to new trajectory models that, when optimized, result in 
significant improvments on the optimality and feasibility over the original Sims-Flanagan method.

\subsection{Continuous thrust time-space propagation}
In this model the impulses at each segment are replaced with continouous thrusts ($F_x$,$F_y$,$F_z$) 
which are assumed to be constant within the segment (see the middle scheme in Figure~\ref{fig:compare3models}). 
Each leg of the trajectory is 
propagated forward and backward with equal-duration segments as before. The propagation of the trajectory changes 
from pure Keplerian to integration of the ordinary differential equations:

$$
\begin{array}{l}
r  =  \sqrt{r_x^2 + r_y^2+r_z^2} \\
F  =  \sqrt{F_x^2 + F_y^2+F_z^2} \\
\dot r_x = v_x\\
\dot r_y = v_y\\
\dot r_z = v_z\\
\dot v_x = - \mu r_x / r^3 + F_x/m \\
\dot v_y = - \mu r_y / r^3 + F_y/m \\
\dot v_z = - \mu r_z / r^3 + F_z/m \\
\dot m = -F/(g_0 I_{sp}) \\
\end{array}
$$

\subsection{Continuous thrust $s$-space propagation}
For the \textsl{continuous thrust $s$-space} method, we apply the Sundman transformation \cite{sundman19112} to change
the independent variable from time $t$ to $s$,and the differential equations in $s$-space becomes:

$$
\begin{array}{l}
r  =  \sqrt{r_x^2 + r_y^2+r_z^2} \\
F  =  \sqrt{F_x^2 + F_y^2+F_z^2} \\
\dot r_x = r v_x\\
\dot r_y = r v_y\\
\dot r_z = r v_z\\
\dot v_x = - \mu r_x / r^2 + r F_x/m \\
\dot v_y = - \mu r_y / r^2 + r F_y/m \\
\dot v_z = - \mu r_z / r^2 + r F_z/m \\
\dot m = -r F/(g_0 I_{sp}) \\
\dot t = r \\
\end{array}
$$

where the derivatives are here meant to be taken with respect to the independent variable $s$. Here each leg of the
trajectory is propagated forward from $s_0$ and backward from $s_f$ in equal $s$-space ($\Delta s$). The time
between the mesh is no longer constant and it is proportional to the radial distance $r$, which implies a
shorter time mesh (a finer grid size or segment) is used when the spacecraft is closer to the central body.
The 8th differential equation gives the condition for matching the time difference between the two endpoints:

\begin{equation}
     T_f-T_0 = \int_{s_0}^{s_f} \! r \, ds
\label{eq:DsConstraint}
\end{equation}  

To implement the \textit{$s$-space} method for optimisation, we assume $s_0$ to be zero and solve for $s_f$ to satisfy Eq.~\ref{eq:DsConstraint}, which means an additional variable and constraint are added for each leg. 

\begin{figure*}[ht]
\centering
\subfloat[][Trajectory plot of an Earth-Mercury rendezvous mission optimised using impulsive $\Delta V$ transcription
	(number of segments = 30).]{
	\includegraphics[width=7.5cm]{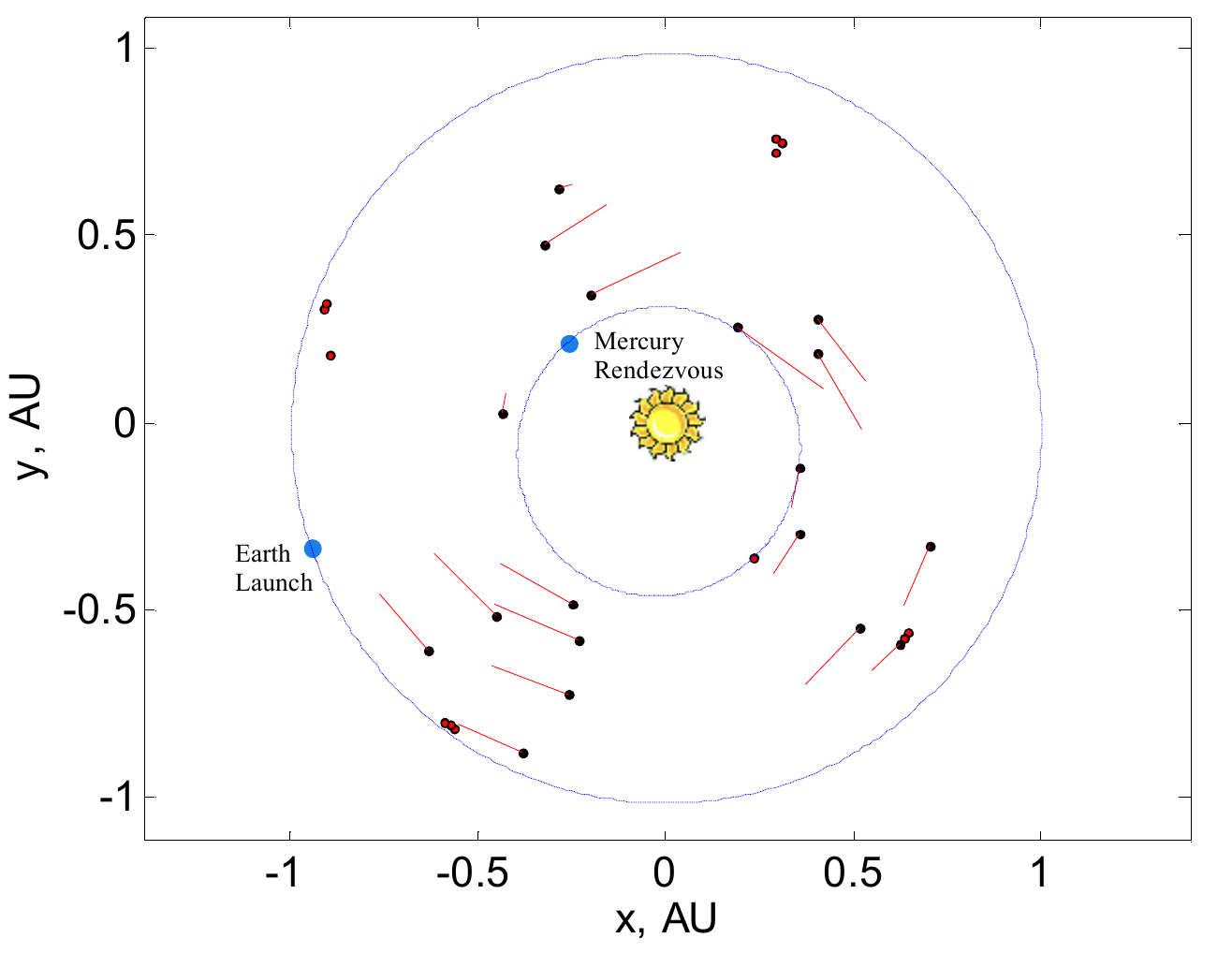}
	\label{fig:EM_IM}
} \quad
\subfloat[][Earth-Mercury trajectory found by continuous thrust time-space propagation.] {
	\includegraphics[width=7.5cm]{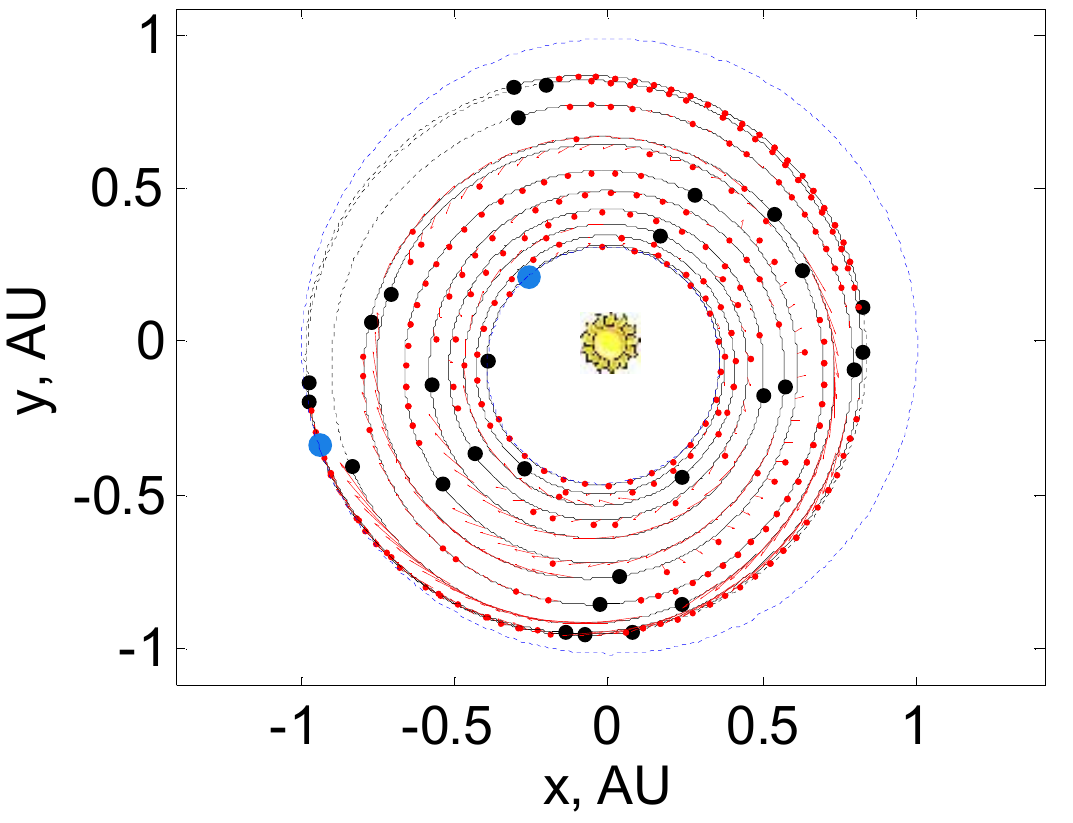}\label{fig:EM_CT}
	}
	\caption{Impulsive and continuous thrust optimal trajectories with equally spaced segments in time.}
\end{figure*}

\begin{figure*}[ht]
\centering
\subfloat[][Earth-Mercury trajectory found by continuous thrust \textit{s}-space propagation.]{
	\includegraphics[width=7.5cm]{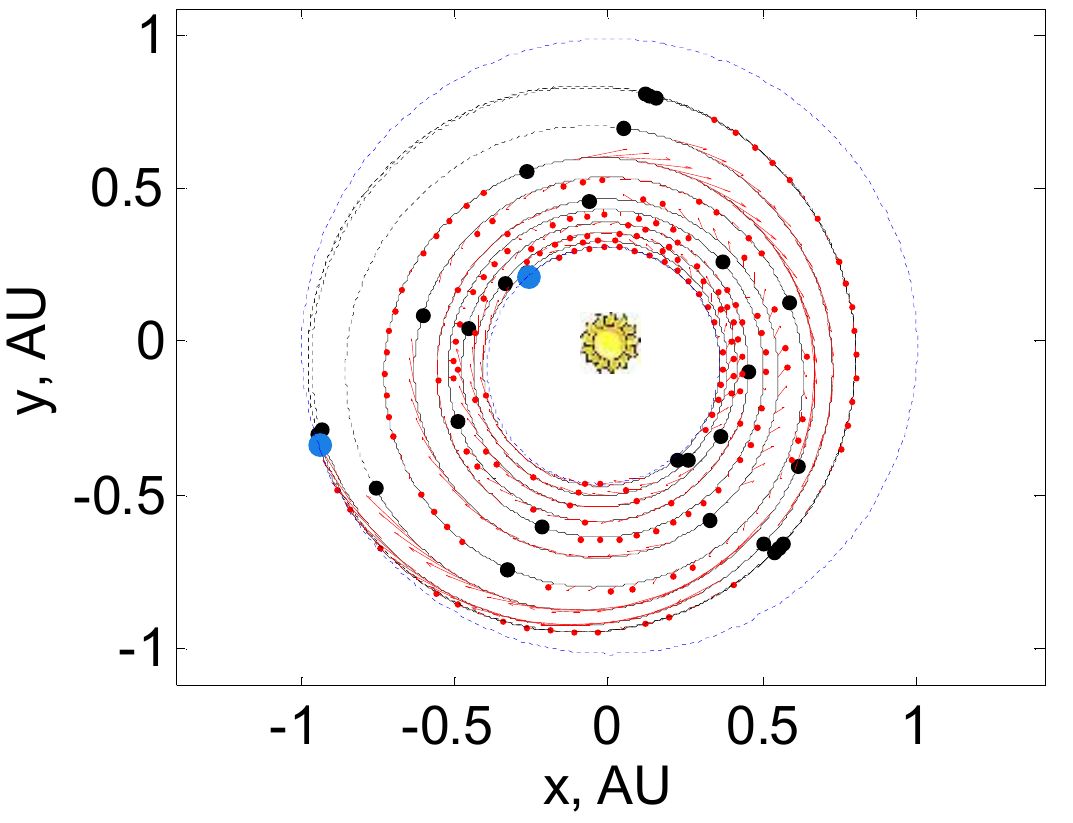}
  \label{fig:EM_Ds}
}\quad
\subfloat[][Duration of eacah segment for the trajectory found by \textit{s}-space propagation.
 The on-line mesh adaptation distributes the segment time duration in a convenient way using the
  Sundman transformation] {
  \includegraphics[width=7.5cm]{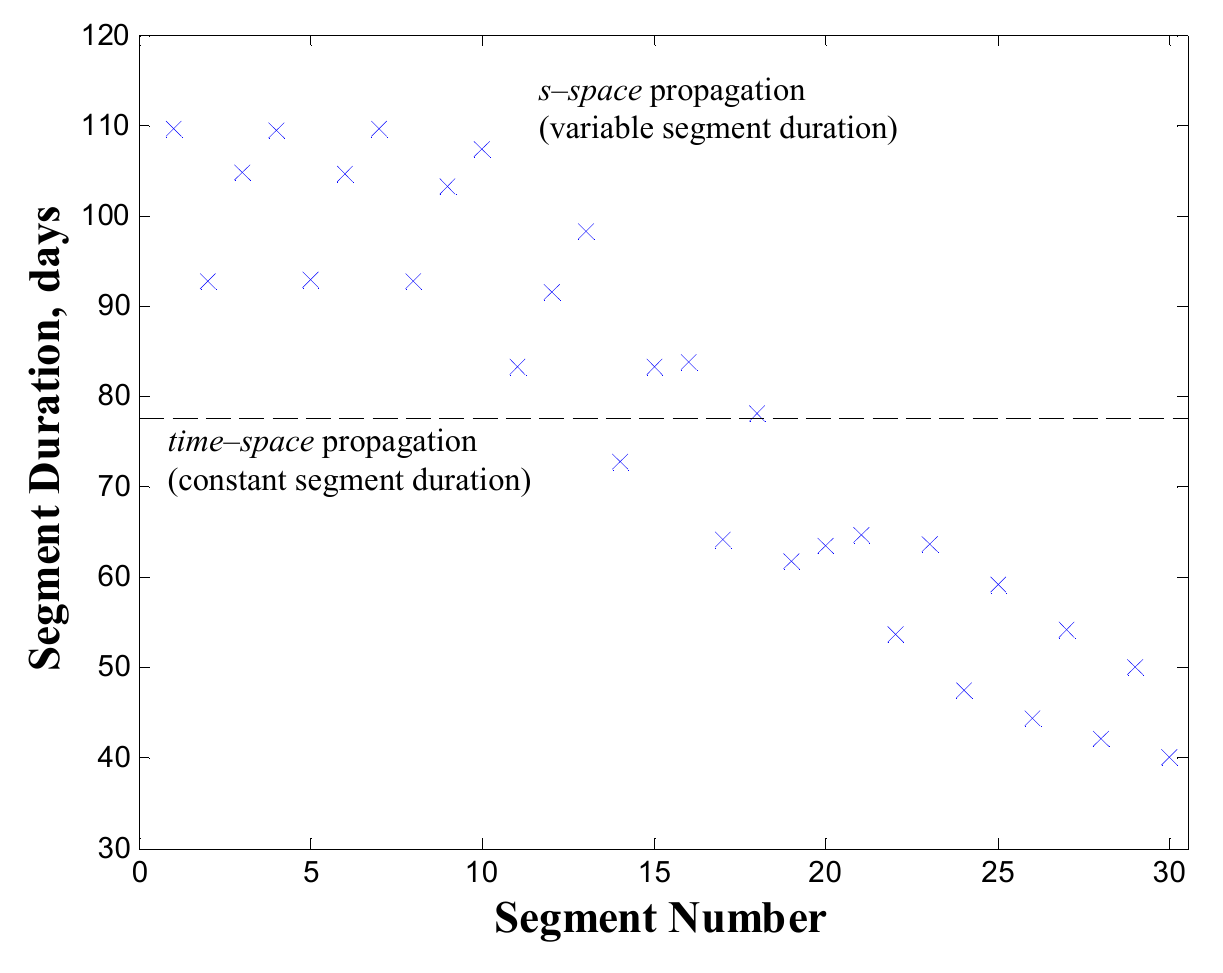}
  \label{fig:SegdurDs}
}
\caption{Trajectory found using the Sundman transformation.}
\end{figure*}

\section{Numerical Example}

We demonstrate our new methods with a sample low-thrust mission to Mercury.  Table~\ref{tab-NEPparam} summarizes the mission specification for a spacecraft similar to the Deep Space 1 mission \cite{Rayman2000}.  Instead of using a model of the SEP (solar electric propulsion) system, we simplify the problem a bit here by assuming a constant thrust and constant specific impulse engine. The launch and arrival dates are kept frozen for the test.

Figure~\ref{fig:EM_IM} shows the trajectory found by the impulsive model, where the black dots
denote the midpoint of the segments and the red lines represent the impulses. It is visually
clear that for this fast rotating trajectory, the impulsive propagation method might not be able
to have a fair representation of the actual low-thrust trajectory. With the same number of segments
(30), the \textsl{continuous thrust time-space model} is able to fill up the "gap" between the
impulses with low-thrust arcs (shown as red curves in figure~\ref{fig:EM_CT}). Even with a lower
final mass (see Table~\ref{tab-EMercuryMf}), the trajectory found by the continuous thrust method
satisfies the real dynamical model which can be used as a trajectory for the actual mission (after
adding other perturbative forces). The trajectory modelling can be further improved by converting to the \textsl{$s$-space} propagation method in figure~\ref{fig:EM_Ds}. In the \textsl{$s$-space} method, we note that when the spacecraft is near Mercury, the grid size is smaller than it is near Earth's orbit, while in the \textsl{time-space} method the segment duration is constant (see figure~\ref{fig:SegdurDs}).
In this example, the implementation of the \textsl{$s$-space} propagation can automatically adapt the
grid size (segment duration) to fit the radial distance (and hence the speed of the spacecraft) of 
the trajectory during the optimisation. A smarter choice of the segment duration along the trajectory
allows the spacecraft to update its control (i.e. thrust) and therefore be more efficient, which explains
why the final mass of the \textsl{$s$-space} is higher than the \textsl{time-space} method.

%



\section{Conclusions and Future Work}

We have successfully extended the Sims-Flanagan model to include the full dynamics of low-thrust trajectory.
The change of independent variable via Sundman transformation can further improve the results through online
adaptive time-mesh during the optimisation. In the future, we hope to investigate a more general form of the
Sundman transformation \cite{Matthew02} and to perform some benchmarking of the new methods to compare
their convergence speed and the accuracy.


\bibliographystyle{abbrv}
\bibliography{sundman_bib}



%
%

\end{document}